\documentclass[leqno]{amsart}
\usepackage{hyperref}
\usepackage{amssymb}
\usepackage{graphicx}
\usepackage{eqnarray,amsthm, amsmath}
\usepackage{mathtools}

\begin{document}
\title[]
{Positive periodic solutions \\ for periodic predator - prey systems \\
of Leslie - Gower or Holling - Tanner type}

\author[]{Dimiter Tsvetkov, Ralitsa Angelova - Slavova}

\address{Dimiter Petkov Tsvetkov (corresponding author) \newline
Department of Communication Systems and Technologies,
Land Forces Faculty,
Vasil Levski National Military University at Veliko Tarnovo,
76 Bulgaria Blvd. \newline
5000 Veliko Tarnovo, Bulgaria}
\email{dimiter.p.tsvetkov@gmail.com}

\address{Ralitsa Lyubomirova Angelova-Slavova \newline
Department of Communication Systems and Technologies,
Land Forces Faculty,
Vasil Levski National Military University at Veliko Tarnovo,
76 Bulgaria Blvd. \newline
5000 Veliko Tarnovo, Bulgaria}
\email{r.angelova.slavova@gmail.com}

%\dedicatory{Communicated by ...}
%\thanks{Submitted .... Published ....}


\begin{abstract}
 In this paper, we consider periodic predator - prey systems of Leslie - Gower or Holling - Tanner type,
 assuming that the coefficients are continuous positive  $\omega$-periodic functions. We prove an existence of positive
 $\omega$-periodic solutions by means of operator method in Banach spaces with cones.
 No constructive conditions are required for the coefficients, besides the positivity and periodicity.
 At the end we point out that the used approach can be applied with minor changes for
proving the existence of positive $\omega$-periodic solutions even in the case of certain delays.
\end{abstract}

\maketitle

\numberwithin{equation}{section}
\newtheorem{theorem}{Theorem}[section]
\newtheorem{prop}{Proposition}[section]
\allowdisplaybreaks

\section{Introduction}
One of the most popular models in mathematical biology is the model about two species, where one of them feeds with another.
It is widely accepted to describe the model dynamics with differential equations.
Recently large amount of publications appear, dedicated to the investigation of certain common properties for the solutions.
The common base of these studies is that they are motivated at all by so called Lotka - Volterra (predator - prey) system,
which appeared first in the prominent monographs of Alfred Lotka \cite{FUND26} and Vito Volterra \cite{FUND37}.

For the general philosophy of the predator - prey systems, besides the mentioned monographs, one can see for example
 Ahmad \& Stamova Eds. \cite{FUND2}, Allen \cite{FUND3}, Bacaer \cite{FUND6},
 Brauer \& Castillo - Chavez \cite{FUNDN2}, Britton \cite{FUND10}, Gillman \cite{FUND15},
 Hadeler \cite{FUND17}, Jones, Plank \& Sleeman \cite{FUND20}, Murray \cite{FUND28}, Smith \cite{FUND31}. Another popular term
 in the related literature is "Kolmogorov system", which applies mostly in the case of several species dynamics, and appears
 to be interchangeable with "Lotka - Volterra system". These systems usually are convenient for being an object of various
 mathematical concepts, from the classic Lyapunov functions, through some topological methods
 (see e.g. Alvarez \& Lazer \cite{TOP5}, Gatica \& Smith \cite{TOP14}, Krasnoselskii \cite{TOP21}, \cite{TOP22},
  Krasnoselskii \& Zabreiko \cite{TOP23}, Mawhin \cite{TOP27}), up to the theory of Lie-Hamilton systems
  (see e.g. Ballesteros at al. \cite{TOP7}, Blasco at al. \cite{TOP9}, de Lucas \& Sardon \cite{TOP13}).

Traditional predator - prey systems are autonomous, including some constant coefficients (parameters),
which have biological meaning (see e.g. Almanza-Vasquez \& Ortiz-Ortiz \cite{AUT4}, Hsu \& Huang \cite{AUT19},
Leslie \& Gower \cite{AUT24}, Tanner \cite{AUT33}, Zhao \cite{AUT40}).
It is reasonable enough to suggest that these coefficients may be replaced by periodic functions of the same period $\omega>0$,
and to look for positive $\omega$-periodic solutions. In this periodic topic, we can mention the works of Ahmad \cite{PER1},
Battauz \& Zanolin \cite{PER8}, Cushing \cite{PER11}, \cite{PER12}, Gopalsamy \cite{PER16}, Pinghua \& Rui \cite{PER29}, Redheffer
\cite{PER30}, Smith \cite{PER32}, Tineo \cite{PER34}, Tineo \& Alvarez \cite{PER35}, Yan \cite{PER38}, Zhang \& Wang \cite{PER39}.
Investigation tools however also change. Often the situation enables to include the operator methods of functional analysis
(see e.g. Lois-Prados \& Precup \cite{PER25}, Tsvetkov \cite{PER36}). Looks like that Krasnoselskii has been the first
who used systematically theory of topological degree to investigate positive solutions of the operator equations
(see \cite{TOP22}, \cite{TOP23}). Note especially his earlier work for cone-compressing or cone-extending operators  \cite{TOP21}.
On this topic one can see also the monographs of Deimling \cite{TOPN3} and Zeidler \cite{TOPN4}.

The mentioned above set of works is not (of course) comprehensive. The reader can discover himself another works,
related to the considered problem (some of them are on open access).

In the present work we use in a very crucial way Krasnoselskii - like result, which is improved essentially
by Gatica \& Smith in \cite{TOP14}.

Hereafter we consider three systems of predator - prey type, as Hsu \& Huang
present them in \cite{AUT19}. For the sake of convenience, we replace the Latin letter notations of the coefficients with Greek ones.
Our attention will be paid to the system of differential equations
\begin{eqnarray*}
\dot{x} &=& \rho x \left(1 - \frac{x}{\kappa} \right) - y p(x) \\
\dot{y} &=& \sigma y \left( 1 - \frac{\eta y}{x} \right)
\end{eqnarray*}
where dot above stands for a differentiation with respect to the independent variable $t$.
In this model the predator (amount $y$) consumes the prey (amount $x$), according to the functional response $p(x)$, where $p(x)$ is a function which takes several shapes.
Consider the following three systems, described according to \cite{AUT19}.

	\textbf{S1}) The functional response is of type 1 with $p(x)=\mu x$. Then we have the following Leslie - Gower model (see also  Leslie \& Gower \cite{AUT24}, Smith \cite{FUND31}).

	\textbf{S2}) The functional response is of type 2 with $p(x) = \frac{\mu x}{\alpha+x}$.
Then we have the following Holling - Tanner model (see also Tanner \cite{AUT33}, Murray \cite{FUND28}, Zhang \cite{AUTN1}).

	\textbf{S3}) The functional response is of type 3 with $p(x) = \frac{\mu x^2}{(\alpha+x)(\beta+x)}$.

	In this work, we shall prove an existence of positive $\omega$-periodic solutions for \textbf{S1}, \textbf{S2} and \textbf{S3} by means of operator method in Banach spaces with cones.

\section{Preliminaries}
For a fixed period $\omega>0$, denote by $\mathcal{C}(\omega) $ the Banach space of continuous $\omega$-periodic functions $f:\mathbb{R} \rightarrow \mathbb{R}$, prowided with the conventional norm
\begin{equation*}
%||f||=||f||_{\infty} = \max_{t}{|f(t)|}
||f||=\max_{t}{|f(t)|}
\end{equation*}
A cone $\mathcal{K}\subset \mathcal{B} $ in a Banach space $\mathcal{B}$ (see e.g. \cite{TOPN3}, \cite{TOP14}, \cite{TOP21},
\cite{TOP22},  \cite{TOP23}, \cite{TOPN4}) is a closed convex set for which $\lambda \mathcal{K} \subset \mathcal{K}$ for $\lambda \ge 0$ and $\mathcal{K}\cap (-\mathcal{K})=\{0\}$ (the zero element of $\mathcal{B}$).
We use the cone $\mathcal{K}(\omega)$ of the nonnegative $\omega$-periodic functions
\begin{equation*}
\mathcal{K}(\omega)=\{f\in \mathcal{C}(\omega):f(t)\ge 0,t \in \mathbb{R}  \}
\end{equation*}
We shall also use cones $\mathcal{K}_\gamma (\omega)\subset\mathcal{K}(\omega)$, $0 <\gamma<1$, of a type
\begin{equation}
 \mathcal{K}_\gamma (\omega) = \{ f\in \mathcal{K}(\omega)|\min_{t}{f(t)} \ge \gamma \max_{t}{f(t)} \}
  \label{eq1}
\end{equation}

%The using of such min-max cone is a keystone in the proof of the main result.

To transform our problem in an operator form, we lay on the following popular proposition, formally borrowed
here from \cite{PER25}.
%Proposition 1
 \begin{prop} \label{Lem1}
	Assume that $a\in \mathcal{C}(\omega), \int_{0}^{\omega}a(\tau)d\tau>0, f\in \mathcal{C} (\omega)$. Then the equation
	\begin{equation*}
	\dot{x} = - a(t)x+f(t)
	\end{equation*}
has unique $\omega$-periodic solutions, which can be obtained in the form
	\begin{equation*}
	x(t)=\int\limits_{t}^{t+\omega}{H(t,s;a)f(s)}ds
    %, t \in \mathbb{R}
	\end{equation*}
where
\begin{equation*}
H(t,s;a)=\frac{ exp(\int_{t}^{s} a(\tau)d\tau)} { exp(\int_{0}^{\omega}a(\tau)d\tau)-1}
\end{equation*}
\end{prop}
The proof is straightforward. $\blacksquare$

At this point we manipulate with a convenient approach, used for example by Lois-Prados \& Precup in \cite{PER25}, to express the solution in an integral operator form. For another integral operator approach one can see Tsvetkov \cite{PER36}.

The role of cone $\mathcal{K}_\gamma (\omega)$ becomes clear from the following proposition, assembled essentially from \cite{TOP22} (see also \cite{PER25}, \cite{PER36}).
%Proposition 2
  \begin{prop} \label{Lem2}
	Assume that $a\in \mathcal{C}(\omega), \int_{0}^{\omega}a(\tau)d\tau>0, f\in \mathcal{K} (\omega)$ and
	\begin{equation}\label{eq0}
	x(t) = \int\limits_{t}^{t+\omega} H(t,s;a)f(s)ds
	\end{equation}
where $H(t,s;a)$ is defined in Proposition \ref{Lem1}. Then $x\in \mathcal{K} (\omega)$. Moreover
	\begin{equation*}
	\min_{t} x(t) \ge \gamma_a \max_{t} x(t)
	\end{equation*}
whit a constant
	\begin{equation*}
	\gamma_a = \frac{\displaystyle \min_{\Omega} H(t,s;a)}{\displaystyle \max_{\Omega} H(t,s;a)}
	\end{equation*}
	where $\Omega = \{ (t,s): 0\le t \le \omega, t \le s \le t+\omega \}$.
Therefore the integral operator, defined in the right side of (\ref{eq0}), maps the cone $\mathcal{K} (\omega)$ into
the cone $\mathcal{K}_{\gamma_a} (\omega)$ (see (\ref{eq1})),
remaining in this way cone $\mathcal{K}_{\gamma_a} (\omega)$ invariant. $\blacksquare$
\end{prop}

Suppose we are given a Banach space $\mathcal{B}$ with a norm $||\cdot||$ and a cone $\mathcal{K}\subset \mathcal{B}$. An operator $\mathcal{T}:D \rightarrow \mathcal{B}, D\subseteq \mathcal{B}$, is said to be compact when $\mathcal{T}$ maps bounded subsets of $D$ into relatively compact subsets of $\mathcal{B}$. Further we shall employ the following basic theorem from Gatica \& Smith \cite{TOP14}.
%Theorem 1.
 \begin{theorem}\label{Th1}
 	(Theorem $1.2$ \cite{TOP14}). Let $ 0 < r < R $ be real numbers
 	\begin{equation*}
 	\mathcal{D} = \{ x \in \mathcal{K}:r \le ||x|| \le R \}
 	\end{equation*}
 	and let $\mathcal{T}:\mathcal{D} \rightarrow \mathcal{K}$ be a compact continuous operator such that
    \begin{align}
    &x\in \mathcal{D}, ||x||=R, \mathcal{T}x=\lambda x \Rightarrow \lambda\le 1& \label{eq:i}\tag{i} \\
    &x\in \mathcal{D}, ||x||=r, \mathcal{T}x=\lambda x \Rightarrow \lambda\ge 1& \label{eq:ii}\tag{ii}\\
    &\inf_{||x||=r} ||\mathcal{T}x|| > 0 \label{eq:iii}\tag{iii}
    \end{align}
    Then $\mathcal{T}$ has a fixed point in $\mathcal{D}$. $\blacksquare$ 	
 	\end{theorem}

Note particularly that in Theorem \ref{Th1}, the cone zero of $\mathcal{B}$ does not belong to the operator domain $\mathcal{D}$.
Note also that for any $f\in\mathcal{K}_\gamma (\omega)$ with $||f||>0$ it holds $f(t)>0$, $\forall t$, i.e.
any nonzero element of $\mathcal{K}_\gamma (\omega)$ is bounded away from zero.

\section{Existence of positive $\omega$-periodic solutions}
In this section we prove that systems \textbf{S1}, \textbf{S2} and \textbf{S3} have positive $\omega$-periodic solutions. Remember that no constructive conditions are required for the coefficients (besides the very natural conditions for the positivity and periodicity). For the sake of clarity we formulate three separate theorems for them in spite of the fact that they are of a quite similar content and thus allow to be gathered in one result.

First pay attention on the system \textbf{S2}, for which our result is the following.
 %Theorem 2.
 \begin{theorem}\label{Th2}
 	Assume that all the coefficients $\rho$, $\kappa$, $\mu$, $\alpha$, $\sigma$ and
 $\eta$ of the system \textbf{S2} (Holling-Tanner model)
 	\begin{eqnarray*}
 		\dot{x} &=& \rho(t) x \left(1 - \frac{x}{\kappa(t)} \right) - \frac{\mu(t) x}{\alpha(t)+x} y \\
 		\dot{y} &=& \sigma(t) y \left( 1 - \frac{\eta(t) y}{x} \right)
 	\end{eqnarray*}	
 are nonnegative continuous $\omega$-periodic functions and non of them is equal to zero identically.
 Assume also that $\alpha(t)\kappa(t) > 0$, $\forall t$, and $\sigma(t)\eta(t)\not\equiv 0$.
 Then system \textbf{S2} has at least one strictly positive $\omega$-periodic solution.
 \end{theorem}

Proof. The proof will be separated in several steps.

\textit{Step 1)} We are looking for strictly positive solutions, therefore we are allowed to change the variables $x=\frac{1}{X}$ and $y=\frac{1}{Y}$, after which the system takes the form
 \begin{eqnarray*}
 	\dot{X} &=& -\rho(t)X +\frac{\rho(t)}{\kappa(t)} + \frac{\mu(t) X^2}{(\alpha(t)X+1)} \frac{1}{Y} \\
 	\dot{Y} &=& -\sigma(t)Y+\sigma(t)\eta(t)X
 \end{eqnarray*}
Then, according to Proposition \ref{Lem1}, the problem for positive $\omega$-periodic solutions of \textbf{S2}
takes the following equivalent operator form
\begin{equation}\label{eq2}
	\begin{aligned}
        &X(t)& &=& &\int\limits_{t}^{t+\omega} H(t,s;\rho) \left(\frac{\rho(s)}{\kappa(s)}+ \frac{\mu(s) X^2(s)}{(\alpha(s)X(s)+1)} \frac{1}{Y(s)}\right)ds \\
       &Y(t)& &=& &\int\limits_{t}^{t+\omega} H(t,s;\sigma) \sigma(s)\eta(s)X(s)ds
	\end{aligned}
\end{equation}
where
\begin{equation*}
H(t,s;\rho)=\frac{exp(\int_{t}^{s} \rho(\tau)d\tau)} {exp(\int_{0}^{\omega}\rho(\tau)d\tau)-1} \quad \text{and} \quad H(t,s;\sigma)=\frac{exp(\int_{t}^{s} \sigma(\tau)d\tau)} {exp(\int_{0}^{\omega}\sigma(\tau)d\tau)-1}
\end{equation*}

Now substitute $Y$ in the first equation of (\ref{eq2}) and obtain a single operator equation
\begin{equation*}
X = \mathcal{T}(X)
\end{equation*}
where the operator $\mathcal{T}(X)$ is given in the right side of (\ref{eq3})
\begin{equation}\label{eq3}
X(t) = \int\limits_{t}^{t+\omega} H(t,s;\rho) \left(\frac{\rho(s)}{\kappa(s)}  + S2Term \right)ds
\end{equation}
with
\begin{equation}\label{S2T}
S2Term = \frac{\mu(s) X^2(s)}{(\alpha(s)X(s)+1)}\frac{1}{\int_{s}^{s+\omega}H(s,\theta;\sigma)\sigma(\theta)\eta(\theta)X(\theta)d\theta}
\end{equation}

Operator $\mathcal{T}(\cdot)$ is defined for all elements of the cone $\mathcal{K}_\gamma (\omega)$, besides the zero, and takes values also in $\mathcal{K}_\gamma (\omega)$, where the cone $\mathcal{K}_\gamma (\omega)$ is defined in (\ref{eq1}) with a constant $\gamma > 0$ which is chosen in the following way (see Proposition \ref{Lem2})
\begin{equation}\label{gam}
\gamma = \frac{\displaystyle\min_{\Omega} H(t,s;\rho)}{\displaystyle\max_{\Omega}H(t,s;\rho)}
\end{equation}

We formulate next proposition manifestly, due to its importance, regardless of its obviousness.
%Proposition 3
\begin{prop} \label{Lem3}
Suppose that $X \in \mathcal{K}_\gamma (\omega)$, $||X|| > 0 $, satisfies the operator equation $X = \mathcal{T}(X)$, given in (\ref{eq3}). Then
\begin{equation*}
Y(t)=\int\limits_{t}^{t+\omega} H(t,\theta;\sigma)\sigma(\theta)\eta(\theta)X(\theta)d\theta
\end{equation*}
and
\begin{equation*}
X(t) =\int\limits_{t}^{t+\omega} H(t,s;\rho) \left(\frac{\rho(s)}{\kappa(s)}+ \frac{\mu(s) X^2(s)}{(\alpha(s)X(s)+1)} \frac{1}{Y(s)}\right)ds
\end{equation*}
satisfy system (\ref{eq2}), therefore $x(t)=\frac{1}{X(t)}$  and $y(t)=\frac{1}{Y(t)}$  are positive $\omega$-periodic solutions of \textbf{S2}. $\blacksquare$
\end{prop}

At this point, our problem reduces to finding a solution of (\ref{eq3}).
We are going to find nonzero  solutions $X \in \mathcal{K}_\gamma (\omega) $ of (\ref{eq3}) by means of Theorem \ref{Th1}.

\textit{Step 2)} Define the basic domain
\begin{equation}\label{bd}
\mathcal{D}=\{ X \in \mathcal{K}_\gamma(\omega):r \le ||X|| \le R\}
\end{equation}
with $\gamma$ from (\ref{gam}). In this case the constants $r$ and $R$ are chosen as follows
\begin{equation}\label{eq4}
0 < r < \max_{t} \int\limits _{t}^{t+\omega} H(t,s;\rho) \frac{\rho(s)}{\kappa(s)}ds
\end{equation}
\begin{equation}\label{eq5}
R > \max_{t} \int\limits _{t}^{t+\omega} H(t,s;\rho) \left( \frac{\rho(s)}{\kappa(s)} +\frac{\mu(s)}{\alpha(s)\gamma^2}   \frac{1}{\int_{s}^{s+\omega}H(s,\theta;\sigma)\sigma(\theta)\eta(\theta)d\theta}\right)ds
\end{equation}

Remember that integrals in (\ref{eq4}) and (\ref{eq5}) are $\omega$-periodic.
Before continue we set yet another trivial proposition, given manifestly due to its importance.
%Proposition 4
\begin{prop} \label{Lem4}
	The operator $\mathcal{T}(X):\mathcal{D} \rightarrow \mathcal{K}_\gamma (\omega)$ is continuous and compact (completely continuous in the terminology of \cite{TOP22}, \cite{TOP23}). $\blacksquare$
\end{prop}

\textit{Step 3.1)} Let us show the validity of clause (\ref{eq:ii}) of Theorem \ref{Th1}. For the sake of contradiction suppose that for some $X \in \mathcal{D}$ with $||X||=r$ and $\mathcal{T}(X)= \lambda X$  we have $\lambda < 1$.
Then $X \ge \lambda X = \mathcal{T}(X)$, which implies
\begin{equation*}
X(t) \ge \int\limits_{t}^{t+\omega} H(t,s;\rho) \left(\frac{\rho(s)}{\kappa(s)} + S2Term
%\frac{\mu(s) X^2(s)}{(\alpha(s)X(s)+1)} \frac{1}{\int_{s}^{s+\omega}H(s,\theta;\sigma)\sigma(\theta)\eta(\theta)X(\theta)d\theta}
\right)ds, \forall t
\end{equation*}
whence, taking into account that $S2Term$ from (\ref{S2T}) is nonnegative, we find
\begin{equation*}
X(t) \ge \int\limits_{t}^{t+\omega} H(t,s;\rho) \frac{\rho(s)}{\kappa(s)}ds (\forall t)
\Rightarrow r=||X|| \ge \max_{t} \int\limits_{t}^{t+\omega} H(t,s;\rho)\frac{\rho(s)}{\kappa(s)}ds
\end{equation*}
that contradicts to the choice of $r$ in (\ref{eq4}).

\textit{Step 3.2)} Now show the validity of clause (\ref{eq:i}) of Theorem \ref{Th1}. For the sake of contradiction suppose that for some $X \in \mathcal{D}$ with $||X||=R$ and $\mathcal{T}(X)=\lambda X$  we have $\lambda > 1$.
Then $X \le \lambda X = \mathcal{T}(X)$, which implies
\begin{equation}\label{eq6}
X(t) \le \int\limits_{t}^{t+\omega} H(t,s;\rho) \left(\frac{\rho(s)}{\kappa(s)} + S2Term
\right)ds, \forall t
\end{equation}
\indent We have $\displaystyle \min_{t}X(t) \ge \gamma R$. Therefore
\begin{equation*}
\begin{split}
S2Term = \frac{\mu(s) X^2(s)}{(\alpha(s)X(s)+1)}
\frac{1}{\int_{s}^{s+\omega}H(s,\theta;\sigma)\sigma(\theta)\eta(\theta)X(\theta) d\theta}\\
\le  \frac{\mu(s) R^2}{(\alpha(s) \gamma R+1)} \frac{1}{\gamma R \int_{s}^{s+\omega}H(s,\theta;\sigma)\sigma(\theta)\eta(\theta)d\theta}\\
\le \frac{\mu(s)}{\alpha(s)\gamma^2} \frac{1}{\int_{s}^{s+\omega}H(s,\theta;\sigma)\sigma(\theta)\eta(\theta)d\theta}
\end{split}
\end{equation*}
Replace in (\ref{eq6}) and obtain
\begin{equation*}\label{eq7}
\begin {split}
R \le \max_{t} \int\limits_{t}^{t+\omega} H(t,s;\rho) \left(\frac{\rho(s)}{\kappa(s)} + S2Term
\right)ds \\
\le \max_{t} \int\limits_{t}^{t+\omega} H(t,s;\rho) \left(\frac{\rho(s)}{\kappa(s)}+ \frac{\mu(s)}{\alpha(s)\gamma^2} \frac{1}{\int_{s}^{s+\omega}H(s,\theta;\sigma)\sigma(\theta)\eta(\theta)d\theta}\right)ds
\end{split}
\end{equation*}
which contradicts to the choice of $R$ in (\ref{eq5}).

\textit{Step 4)} Finally show the validity of clause (\ref{eq:iii}) of Theorem \ref{Th1}. Let $X \in \mathcal{K}_\gamma (\omega)$ and $||X||=r$. Then we get
\begin{equation*}
||\mathcal{T}(X)||=\max_{t}\mathcal{T}(X)(t)
\ge \max_{t} \int\limits_{t}^{t+\omega} H(t,s;\rho) \frac{\rho(s)}{\kappa(s)} ds = Const > 0
\end{equation*}
where $Const$ does not depend on $X$. Therefore
\begin{equation*}
\inf_{||X||=r}||\mathcal{T}X|| \ge Const >0
\end{equation*}
Theorem \ref{Th2} is proved. $\blacksquare$

Continue with system \textbf{S1}.
%Theorem 3.
\begin{theorem}\label{Th3}
	Assume that all the coefficients $\rho$, $\kappa$, $\mu$, $\sigma$ and $\eta$ of the system \textbf{S1} (Leslie-Gower  model)
	\begin{eqnarray*}
		\dot{x} &=& \rho(t) x \left(1 - \frac{x}{\kappa(t)} \right) -  \mu(t) xy \\
		\dot{y} &=& \sigma(t) y \left( 1 - \frac{\eta(t) y}{x} \right)
	\end{eqnarray*}	
are nonnegative continuous $\omega$-periodic functions and non of them is equals to zero identically.
Assume also that $\kappa(t) > 0$, $\forall t$, and $\sigma(t)\eta(t)\not\equiv 0$.
Then system \textbf{S1} has at least one strictly positive $\omega$-periodic solution.
\end{theorem}

Proof. The proof goes in the same way as in Theorem \ref{Th2}. Changing of the variables $x = \frac{1}{X}$ and $y = \frac{1}{Y}$, transforms system \textbf{S1} to
\begin{eqnarray*}
	\dot{X} &=& -\rho(t)X + \frac{\rho(t)}{\kappa(t)}+  \mu(t) \frac{X}{Y} \\
	\dot{Y} &=& -\sigma(t)Y +\sigma(t)\eta(t)X
\end{eqnarray*}	
The corresponding equivalent operator form is
\begin{eqnarray*}
	\dot{X} &=& \int\limits_{t}^{t+\omega}H(t,s;\rho) \left(\frac{\rho(s)}{\kappa(s)}+  \mu(s) \frac{X(s)}{Y(s)} \right)ds\\
	\dot{Y} &=& \int\limits_{t}^{t+\omega}H(t,s;\sigma)\sigma(s)\eta(s)X(s)ds
\end{eqnarray*}
with a corresponding single operator equation on $\mathcal{K}_\gamma (\omega)$
    \begin{equation*}
	X(t) = \int\limits_{t}^{t+\omega}H(t,s;\rho) \left(\frac{\rho(s)}{\kappa(s)} + \mu(s)X(s) \frac{1}{\int\limits_{s}^{s+\omega} H(s,\theta;\sigma)\sigma(\theta)\eta(\theta)X(\theta)d\theta} \right)ds
	\end{equation*}

Consider the basic domain $\mathcal{D}$ from (\ref{bd}),
where the constant $r$ is chosen as in (\ref{eq4}) and $R$ is chosen as follows
\begin{equation*}
R > \max_{t} \int\limits_{t}^{t+\omega}H(t,s;\rho) \left(\frac{\rho(s)}{\kappa(s)}+\mu(s)\frac{1}{\gamma \int_{s}^{s+\omega}H(s,\theta;\sigma)\sigma(\theta)\eta(\theta)d\theta}\right)ds
\end{equation*}

Rest of the proof is similar to the end of the proof of Theorem \ref{Th2}. $\blacksquare$
%Theorem 4.
\begin{theorem}\label{Th4}
	Assume that all the coefficients $\rho$, $\kappa$, $\mu$, $\alpha$, $\beta$, $\sigma$,
and $\eta$ of the system \textbf{S3}
	\begin{eqnarray*}
		\dot{x} &=& \rho(t) x \left(1 - \frac{x}{\kappa(t)} \right) -\frac{\mu(t) x^2}{(\alpha(t)+x)(\beta(t)+x)}y   \\
		\dot{y} &=& \sigma(t) y \left( 1 - \frac{\eta(t) y}{x} \right)
	\end{eqnarray*}	
are nonnegative continuous $\omega$-periodic functions and non of them is equals to zero identically.
Assume also that $\alpha(t)\beta(t)\kappa(t) > 0$, $\forall t$, and $\sigma(t)\eta(t)\not\equiv 0$.
Then system \textbf{S3} has at least one strictly positive $\omega$-periodic solution.
\end{theorem}

Proof. The proof again goes in the same way as in Theorem \ref{Th2}. Changing of the variables $x=\frac{1}{X}$ and $y=\frac{1}{Y}$, transforms system \textbf{S3} to
\begin{eqnarray*}
	\dot{X} &=& -\rho(t)X + \frac{\rho(t)}{\kappa(t)}+ \frac{\mu(t)X^2}{(\alpha(t)X+1)(\beta(t)X+1)}  \frac{1}{Y} \\
	\dot{Y} &=& -\sigma(t)Y +\sigma(t)\eta(t)X
\end{eqnarray*}
The corresponding equivalent operator form is
\begin{eqnarray*}
	X(t) &=& \int\limits_{t}^{t+\omega} H(t,s;\rho)\left( \frac{\rho(s)}{\kappa(s)}+ \frac{\mu(s)X^2(s)}{(\alpha(s)X(s)+1)(\beta(s)X(s)+1)}  \frac{1}{Y(s)} \right)ds\\
	Y(t) &=& \int\limits_{t}^{t+\omega}H(t,s;\sigma)\sigma(s)\eta(s)X(s)ds
\end{eqnarray*}
with a corresponding single operator equation on $\mathcal{K}_\gamma (\omega)$
\begin{equation*}
	X(t) = \int\limits_{t}^{t+\omega} H(t,s;\rho)\left( \frac{\rho(s)}{\kappa(s)} + S3Term \right)ds
\end{equation*}
where
\begin{equation*}
S3Term = \frac{\mu(s)X^2(s)}{(\alpha(s)X(s)+1)(\beta(s)X(s)+1)}  \frac{1}{\int_{s}^{s+\omega}H(s,\theta;\sigma)\sigma(\theta)\eta(\theta)X(\theta)d\theta}
\end{equation*}

Consider again the basic domain $\mathcal{D}$ from (\ref{bd}),
where $r$ is chosen as in (\ref{eq4}) and $R$ is chosen as follows
\begin{equation*}
R> \max_{t} \int\limits_{t}^{t+\omega} H(t,s;\rho)\left(\frac{\rho(s)}{\kappa(s)} + \frac{\mu(s)}{R\alpha(s)\beta(s)\gamma^3} \frac{1}{\int_{s}^{s+\omega}H(s,\theta;\sigma)\sigma(\theta)\eta(\theta)d\theta}\right)ds
\end{equation*}

Rest of the proof again is similar to the end of the proof of Theorem \ref{Th2}. $\blacksquare$

\section{Conclusions and numerical examples}

Using of the operator method and cone technique for finding positive solutions has a history, starting from the works of Krasnoselskii (see e.g. \cite{TOP21}, \cite{TOP22}, \cite{TOP23}).
%Using of the min-max cone $\mathcal{K}_\gamma (\omega)$ is the most important detail in the proof.
The using of min-max cone $\mathcal{K}_\gamma (\omega)$ is the keystone in the proof of the main result.
It helps to control the growth of the resolving operator.
Cones of this type appear for example in \cite{TOP22}, \cite{TOP23}.

This min-max cone approach is applied by Tsvetkov in \cite{PER36} for the very classical Lotka - Volterra predator - prey model
\begin{eqnarray*}
	\dot{x} &=&  \alpha x - \lambda xy\\
	\dot{y} &=& -\beta y + \chi xy
\end{eqnarray*}
where the coefficients $\alpha$, $\beta$, $\lambda$ and $\chi$ are nonnegative $\omega$-periodic functions,
none of which is equal to zero identically.

The proof scheme in \cite{PER36}, however does not fit the predator - prey model with logistic growth on prey
\begin{eqnarray*} %\label{pplg}
	\dot{x} &=& \alpha x\left( 1-\frac{x}{\kappa}\right) - \lambda xy\\
	\dot{y} &=& -\beta y+\chi xy
\end{eqnarray*}
because the existence technics therein cannot remove the case $y \equiv 0$ in the latter system.

The min-max cone approach is also applied by Prados \& Precup \cite{PER25} for Lotka - Volterra systems considered therein.

Our method includes the autonomous systems in a trivial way.
In this case, the constant coefficients are periodic functions of any period, and the constant solution may serve as
an $\omega$-periodic solution for the problem considered.
One may say that an usual $\omega$-periodic solution is born from the constant solution,
when the coefficients turn from constants to $\omega$-periodic functions. In the autonomous case, any of
our three systems \textbf{S1}, \textbf{S2} and \textbf{S3} has a positive constant solution.

Keep in mind that an autonomous Lotka - Volterra system may exhibit very remarkable behavior (see e.g. Hofbauer \& So \cite{AUT18}).
We finish this work with three illustrative numerical examples.

% Example 1.
\textbf{Example 1.} Consider the following system of type \textbf{S2}
\begin{eqnarray*}
	\dot{x} &=& (1+\sin{5t})x\left(1-\frac{x}{2+\sin{t}}\right)-\frac{(1+\cos{3t})x}{(2-\cos{3t})+x}y\\
	\dot{y} &=& (1-\cos{7t})y\left(1-\frac{(1-\sin{t})y}{x}\right)
\end{eqnarray*}
A $2\pi$-periodic solution is found near the initial values
\begin{equation*}
	x(0) = 0.8416874693971644, y(0) = 0.5259233975099778
\end{equation*}
with
\begin{equation*}
	|x(0)-x(2\pi)| + |y(0)-y(2\pi)| < 1.5096 \times 10^{-11}
\end{equation*}
\begin{figure}[h]
\includegraphics[scale=0.2,bb=0 0 1600 400]{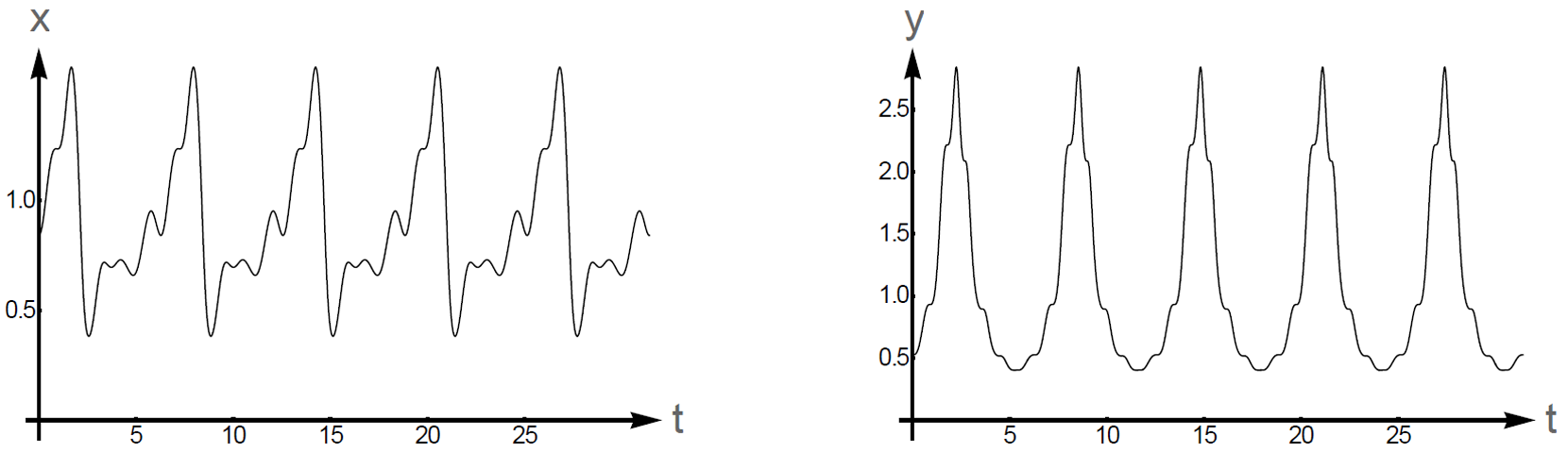}
\caption{Positive $2\pi$-periodic solution for Example 1 for $t\in[0,10\pi]$. }
\end{figure}

%Example 2.
\textbf{Example 2.} Consider the following system of type \textbf{S3}
\begin{eqnarray*}
	\dot{x} &=& (1+\sin{t})x\left(1-\frac{x}{2+\sin{2t}}\right)-\frac{(2+\cos{3t})x^2}{(2+\cos{2t}+x)(2-\cos{3t}+x)}y\\
	\dot{y} &=& (1+\cos{t})y\left(1-\frac{(1-\sin{t})y}{x}\right)
\end{eqnarray*}
A $2\pi$-periodic solution is found near the initial values
\begin{equation*}
	x(0) = 0.6406510789582541, y(0) = 0.4091984714503302
\end{equation*}
with
\begin{equation*}
	|x(0)-x(2\pi)| + |y(0)-y(2\pi)| < 3.1943 \times 10^{-11}
\end{equation*}
\begin{figure}[h]
\includegraphics[scale=0.2,bb=0 0 1600 400]{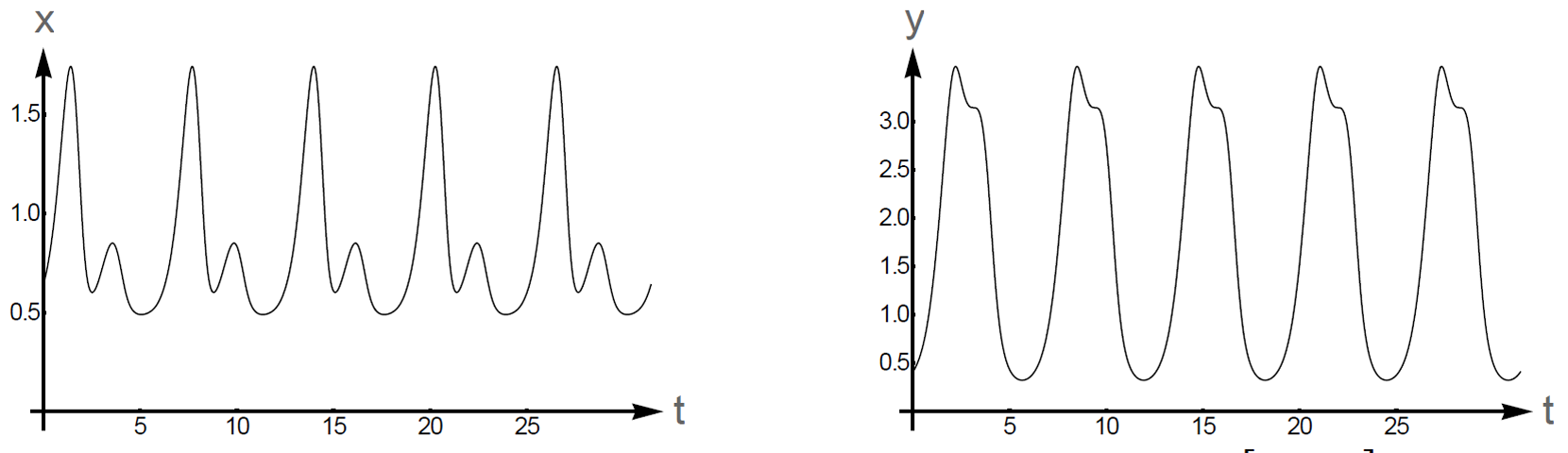}
\caption{Positive $2\pi$-periodic solution for Example 2 for $t\in[0,10\pi]$. }
\end{figure}

%Example 3.
\textbf{Example 3.} Consider the following system of type \textbf{S1}
\begin{eqnarray*}
	\dot{x} &=& (1+\sin{2t})x\left(1-\frac{x}{2+\sin{5t}}\right)-(1+\cos{3t})xy\\
	\dot{y} &=& (1-\cos{t})y\left(1-\frac{(1-\sin{t})y}{x}\right)
\end{eqnarray*}
A $2\pi$-periodic solution is found near the initial values
\begin{equation*}
	x(0) = 0.6504022496685088, y(0) = 0.3825388660004428
\end{equation*}
with
\begin{equation*}
	|x(0)-x(2\pi)| + |y(0)-y(2\pi)| < 3.3083 \times 10^{-11}
\end{equation*}

% Figure 3.
\begin{figure}[h]
\includegraphics[scale=0.2,bb=0 0 1600 400]{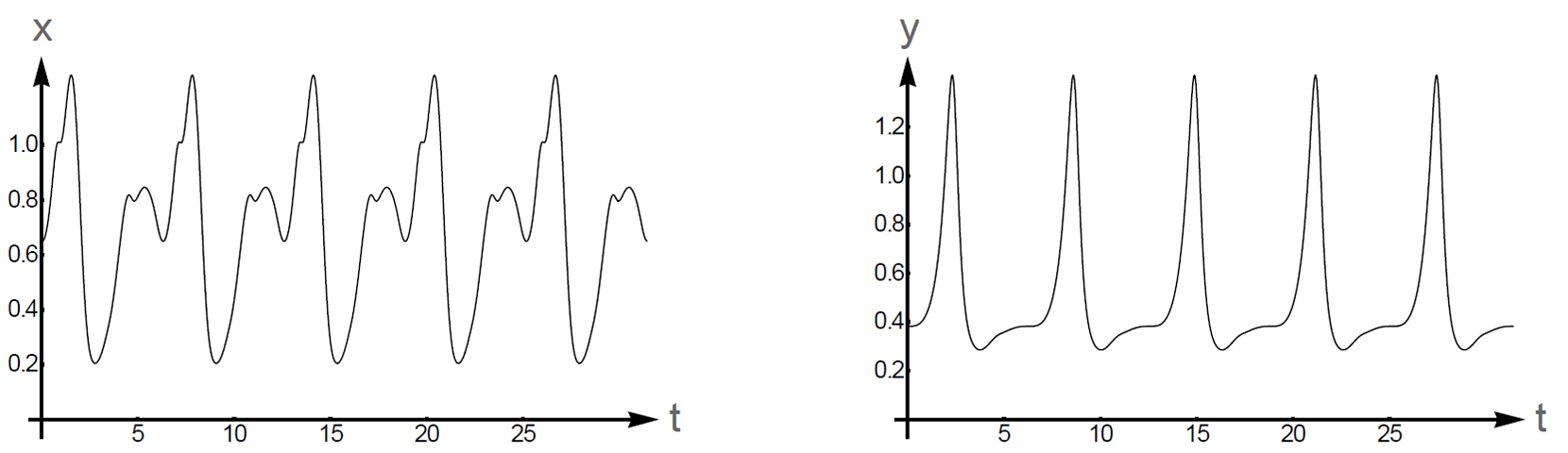}
\caption{Positive $2\pi$-periodic solution for Example 3 for $t\in[0,10\pi]$. }
\end{figure}

In all three examples above the solution appears locally stable, but our operator method do no allow itself
to investigate the stability. On the other hand, the operator method is very convenient for
proving the existence of positive periodic solutions even in the case of delays.
For instance, one can find that the conditions of Theorem \ref{Th2}
give existence of positive $\omega$-periodic solutions of the following delay (Holling-Tanner) system
	\begin{eqnarray*}
 		\dot{x}(t) &=& \rho(t) x(t) \left(1 - \frac{x(t)}{\kappa(t)} \right) - \frac{\mu(t) x(t)}{\alpha(t)+x(t)} y(t-\tau_y) \\
 		\dot{y}(t) &=& \sigma(t) y(t) \left( 1 - \frac{\eta(t) y(t)}{x(t-\tau_x)} \right)
 	\end{eqnarray*}	
with some constants $\tau_x\geq 0$ and $\tau_y\geq 0$. In this case $S2Term$ from (\ref{S2T}) takes the form
\begin{equation*}
\frac{\mu(s) X^2(s)}{(\alpha(s)X(s)+1)}\frac{1}{\int_{s-\tau_y}^{s-\tau_y+\omega}H(s-\tau_y,\theta;\sigma)
\sigma(\theta)\eta(\theta)X(\theta-\tau_x)d\theta}
\end{equation*}
which in fact implies no substantially alterations in the proof.

%\subsection*{Acknowledgements}


\begin{thebibliography}{99}
\bibitem{PER1} S. Ahmad;
\emph{Convergence and Ultimate Bounds of Solutions of the Nonautonomous Volterra - Lotka Competition Equations},
J. Math. Anal. Appl., \textbf{127} (1987), 377 - 387.

\bibitem{FUND2} S. Ahmad, I. Stamova (Eds.);
 \emph{Lotka-Volterra and Related Systems: Recent Developments in Population Dynamics},
De Gruyter, Berlin/Boston, 2013. ISBN 978-3-11-026951-2

\bibitem{FUND3} L.J.S. Allen;
 \emph{An introduction to mathematical biology},
 Upper Saddle River: Pearson education, 2007. ISBN 0-13-035216-0
	
\bibitem{AUT4} E. Almanza - Vasquez, R - D. Ortiz - Ortiz;
\emph{Bifurcations in the Dynamics of Rosenzweig-MacArthur Predator-Prey Model Considering Saturated Refuge for the Preys},
Applied Mathematical Sciences, \textbf{9} (2015), 7475 - 7482.

\bibitem{TOP5} C. Alvarez, A.C. Lazer;
\emph{An application of topological degree to the periodic competing species problem},
J. Austral. Math. Soc. Ser. B, \textbf{28} (1986), 202 - 219.

\bibitem{FUND6} N. Bacaer;
 \emph{A Short History of Mathematical Population Dynamics},
Springer, New York, 2010. ISBN 978-0-85729-114-1

\bibitem{TOP7} A. Ballesteros, A. Blasco, F.J. Herranz, J. de Lucas, C. Sardon;
\emph{Lie-Hamilton systems on the plane: properties, classification and applications},
J. Diff. Eqns., \textbf{258} (2015), no. 8, 2873 - 2907.

\bibitem{PER8} A. Battauz, F. Zanolin;
\emph{Coexistence States for Periodic Competitive Kolmogorov Systems},
J. Math. Ann. Appl., \textbf{219} (1998), 179 - 199.

\bibitem{TOP9} A. Blasco, F.J. Herranz, J. de Lucas, C. Sardon;
\emph{Lie - Hamilton systems on the plane: applications and superposition rules},
J. Phys. A, \textbf{48} (2015), 345 - 202.

\bibitem{FUNDN2} F. Brauer, C. Castillo - Chavez;
\emph{Mathematical models in population biology and epidemiology},
Springer, New-York, 2012. ISBN 978-1-4614-1685-2

\bibitem{FUND10} N. F. Britton;
 \emph{Essential Mathematical Biology},
Springer-Verlag, London, 2003. ISBN 978-1-85233-536-6

\bibitem{PER11} J. Cushing;
\emph{Periodic Kolmogorov Systems},
SIAM J. Math. Anal., \textbf{13} (1982), No 5, 811 - 827.

\bibitem{PER12} J. Cushing;
\emph{Two Species Competition in a Periodic Environment},
J. Math. Biology, \textbf{10} (1980), 385 - 400.

\bibitem{TOPN3} K. Deimling;
 \emph{Nonlinear Functional Analysis},
Springer - Verlag, Berlin Heidelberg, 1985. ISBN 0-387-13928-1
	
\bibitem{TOP13} J. de Lucas, C. Sardon;
 \emph{A Guide To Lie Systems With Compatible Geometric Structures},
New Jersey: World Scientific, 2020. ISBN-10 1786346974

\bibitem{TOP14} J. A. Gatica, H.L. Smith;
\emph{Fixed point techniques in a cone with applications},
J. Math. Anal. Appl., \textbf{61} (1977), 58 - 71.

\bibitem{FUND15} M. Gillman;
 \emph{An Introduction to Mathematical Models in Ecology and Evolution: Time and Space},
Wiley-Blackwell, New York, 2009. ISBN-10 140517515X

\bibitem{PER16} K. Gopalsamy;
\emph{Global Asymptotic stability in a periodic Lotka - Volterra System},
J. Austral. Math. Soc. Ser. B, \textbf{27} (1985), 66 - 72.

\bibitem{FUND17} K. P. Hadeler;
 \emph{Topics in Mathematical Biology},
Springer International Publishing AG, 2017. ISBN 978-3-319-65620-5

\bibitem{AUT18} J. Hofbauer, J. W-H. So;
\emph{Multiple Limit Cycles for Three Dimensional Lotka-Volterra Equations},
Appl. Math. Lett., \textbf{7} (1994), 65 - 70.
	
\bibitem{AUT19} S. Hsu, T. Huang;
\emph{Global Stability for a Class of Predator - Prey Systems},
SIAM J. Appl. Math., \textbf{55} (1995), 763 - 783.
	
\bibitem{FUND20} D. S. Jones, M. J. Plank, B. D. Sleeman;
 \emph{Differential Equations and Mathematical Biology},
Boca Raton, FL: Chapman \& Hall. ISBN 1-4200-8358-9

\bibitem{TOP21} M. A. Krasnoselskii;
\emph{Fixed points of cone-compressing or cone-extending operators},
Dokl. Akad. Nauk SSSR, \textbf{135} (1960), 527 - 530.

\bibitem{TOP22} M. A. Krasnoselskii;
 \emph{Positive solutions of operator equations. Translated from the Russian by Richard E. Flaherty},
Noordhoff, Groningen, 1966.
	
\bibitem{TOP23} M. A. Krasnoselskii, P.P. Zabreiko;
 \emph{Geometrical methods of nonlinear analysis. Translated from the Russian by Christian C. Fenske},
Springer - Verlag, Berlin, 1984. ISBN 3-540-12945-6

\bibitem{AUT24} P. H. Leslie, J.C. Gower;
\emph{The properties of stochastic model for the predator-prey type of interaction between two species},
Biometrika, \textbf{47} (1960), 219 - 234.
	
\bibitem{PER25} C. Lois-Prados, R. Precup;
\emph{Positive periodic solutions for Lotka - Volterra systems with a general attack rate},
Nonlinear Analysis: Real World Applications, \textbf{52} (2020).

\bibitem{FUND26} A. J. Lotka;
 \emph{Elements of Physical Biology},
Williams and Wilkins, Baltimore, 1925.

\bibitem{TOP27} J. Mawhin;
\emph{Leray - Schauder continuation theorems in the absence of a priori bounds},
Topological Methods in Nonlinear Analysis, Journal of the Juliusz Schauder Center, \textbf{9} (1997), 179 - 200.

\bibitem{FUND28} J. D. Murray;
 \emph{Mathematical Biology},
Springer-Verlag, New York, 2002. ISBN 0-387-95223-3

\bibitem{PER29} Y. Pinghua, X. Rui;
\emph{Global Attractivity of the Periodic Lotka - Volterra System},
J. Math. Anal. Appl., \textbf{233} (1999), 221 - 232.

\bibitem{PER30} R. Redheffer;
\emph{Nonautonomous Lotka - Volterra Systems, I},
J. Diff. Eqns., \textbf{127} (1996), 519 - 541.

\bibitem{FUND31} J. M. Smith;
 \emph{Models in ecology},
University Press, Cambridge, 1978. ISBN 0-521-29440-1

\bibitem{PER32} H. L. Smith;
\emph{Periodic Orbits of Competitive and Cooperative Systems},
J. Diff. Equations, \textbf{65} (1986), 361 - 373.

\bibitem{AUT33} J. T. Tanner;
\emph{The stability and the intrinsic growth rates of prey and predator populations},
Ecology, \textbf{56} (1975), 855 - 867.

\bibitem{PER34} A. Tineo;
\emph{An Iterative Scheme for the N-Competing Species Problem},
J. Diff. Eqns., \textbf{116} (1995), 1 - 15.

\bibitem{PER35} A. Tineo, C. Alvarez;
\emph{A different consideration about the globally asymptotically stable solution of the periodic $n$-competing species problem},
J. Math. Anal. Appl., \textbf{159} (1991), 44 - 50.

\bibitem{PER36} D. P. Tsvetkov;
\emph{A Periodic Lotka - Volterra System},
Serdica Math. J., \textbf{22} (1996), 109 - 116.

\bibitem{FUND37} V. Volterra;
 \emph{Theory of functionals and of integral and integro-differential equations. With a preface by G. C. Evans, a biography of Vito Volterra and a bibliography of his published works by E. Whittaker},
Dover Publications, Inc., New York, 1959.

\bibitem{PER38} J. Yan;
\emph{Global positive periodic solutions of periodic $n$-species competition systems},
J. Math. Anal. Appl., \textbf{356} (2009), 288 - 294.

\bibitem{TOPN4} E. Zeidler;
 \emph{Nonlinear Functional Analysis and Applications I. Fixed-Point Theorems},
Springer - Verlag, New York, 1986. ISBN 978-0-387-90914-1

\bibitem{AUTN1} J-F. Zhang;
\emph{Bifurcation analysis of a modified Holling - Tanner predator - prey model with time delay},
Applied Mathematical Modelling, \textbf{36} (2012), 1219 - 1231.

\bibitem{PER39} Z. Zhang, Z. Wang;
\emph{Periodic solutions for nonautonomous predator-prey system with diffusion and time delay},
Hiroshima Math. J., \textbf{31} (2001), 371 - 381.

\bibitem{AUT40} X-Q. Zhao;
\emph{The qualitative analysis of $n$-species Lotka-Volterra periodic competition systems},
Mathematical and Computer Modelling, \textbf{15} (1991), 3 - 8.

\end{thebibliography}
\end{document}